\documentstyle[amssymb,12pt]{amsart}
\newtheorem{prop}{Proposition}[section]
\newtheorem{prop:def}{Proposition-Definition}[section]

\newtheorem{lemma}{Lemma}[section]
\newtheorem{thm}{Theorem}[section]
\newtheorem{cor}{Corollary}[section]

\theoremstyle{remark}
\newtheorem{remark}{Remark}

\begin{document}

\newcommand{\nc}{\newcommand} \nc{\on}{\operatorname}

\nc{\pa}{\partial}

\nc{\cA}{{\cal A}} \nc{\cB}{{\cal B}}\nc{\cC}{{\cal C}} 
\nc{\cD}{{\cal D}} 
\nc{\cE}{{\cal E}} \nc{\cG}{{\cal G}}\nc{\cH}{{\cal H}} 
\nc{\cI}{{\cal I}} \nc{\cJ}{{\cal J}}\nc{\cK}{{\cal K}} 
\nc{\cL}{{\cal L}} \nc{\cR}{{\cal R}} \nc{\cS}{{\cal S}}   
\nc{\cV}{{\cal V}} \nc{\cX}{{\cal X}}

\nc{\sh}{\on{sh}}\nc{\Id}{\on{Id}}\nc{\Diff}{\on{Diff}}
\nc{\Perm}{\on{Perm}}\nc{\conc}{\on{conc}}\nc{\Alt}{\on{Alt}}
\nc{\ad}{\on{ad}}\nc{\Der}{\on{Der}}\nc{\End}{\on{End}}
\nc{\no}{\on{no\ }} \nc{\res}{\on{res}}\nc{\ddiv}{\on{div}}
\nc{\Sh}{\on{Sh}} \nc{\card}{\on{card}}\nc{\dimm}{\on{dim}}
\nc{\Sym}{\on{Sym}} \nc{\Jac}{\on{Jac}}\nc{\Ker}{\on{Ker}}
\nc{\Vect}{\on{Vect}} \nc{\Spec}{\on{Spec}}\nc{\Cl}{\on{Cl}}
\nc{\Imm}{\on{Im}}\nc{\limm}{\lim}\nc{\Ad}{\on{Ad}}
\nc{\ev}{\on{ev}} \nc{\Hol}{\on{Hol}}\nc{\Det}{\on{Det}}
\nc{\Bun}{\on{Bun}}\nc{\diag}{\on{diag}}\nc{\pr}{\on{pr}} 
\nc{\Span}{\on{Span}}\nc{\Comp}{\on{Comp}}\nc{\Part}{\on{Part}}
\nc{\tensor}{\on{tensor}}\nc{\ind}{\on{ind}}\nc{\id}{\on{id}}
\nc{\Hom}{\on{Hom}}\nc{\Quant}{\on{Quant}}\nc{\Def}{\on{Def}}
\nc{\AutLBA}{\on{AutLBA}}\nc{\AutQUE}{\on{AutQUE}}
\nc{\LBA}{\on{LBA}}\nc{\Aut}{\on{Aut}}\nc{\QUE}{\on{QUE}}
\nc{\Lyn}{\on{Lyn}}\nc{\LCA}{\on{LCA}}\nc{\FLBA}{\on{FLBA}}
\nc{\QuantFLBA}{\on{QuantFLBA}}\nc{\QTA}{\on{QTA}}
\nc{\LA}{\on{LA}}\nc{\FLA}{\on{FLA}}\nc{\EK}{\on{EK}}
\nc{\class}{\on{class}}\nc{\br}{\on{br}}\nc{\co}{\on{co}}
\nc{\Prim}{\on{Prim}}\nc{\ren}{\on{ren}}\nc{\lbr}{\on{lbr}}
\nc{\SC}{\on{SC}}

\nc{\al}{\alpha}\nc{\g}{\gamma}\nc{\de}{\delta}
\nc{\eps}{\epsilon}\nc{\la}{{\lambda}}
\nc{\si}{\sigma}\nc{\z}{\zeta}

\nc{\La}{\Lambda}

\nc{\ve}{\varepsilon} \nc{\vp}{\varphi} 

\nc{\AAA}{{\mathbb A}}\nc{\CC}{{\mathbb C}}\nc{\ZZ}{{\mathbb Z}} 
\nc{\QQ}{{\mathbb Q}} \nc{\NN}{{\mathbb N}}\nc{\VV}{{\mathbb V}} 
\nc{\KK}{{\mathbb K}} 

\nc{\ff}{{\mathbf f}}\nc{\bg}{{\mathbf g}}
\nc{\ii}{{\mathbf i}}\nc{\kk}{{\mathbf k}}
\nc{\bl}{{\mathbf l}}\nc{\zz}{{\mathbf z}} 
\nc{\pp}{{\mathbf p}}\nc{\qq}{{\mathbf q}}

\nc{\cF}{{\cal F}}\nc{\cM}{{\cal M}}\nc{\cO}{{\cal O}}
\nc{\cT}{{\cal T}}\nc{\cW}{{\cal W}}
\nc{\Assoc}{{\mathbf Assoc}}

\def\sha{{\mathop{\scriptstyle\sqcup\!\hspace{-0.2pt}\sqcup}}}
\def\Sha{{\mathop{\scriptstyle\amalg\!\hspace{-1.52pt}\amalg}}}


\nc{\ub}{{\underline{b}}}
\nc{\uk}{{\underline{k}}} \nc{\ul}{{\underline}}
\nc{\un}{{\underline{n}}} \nc{\um}{{\underline{m}}}
\nc{\up}{{\underline{p}}}\nc{\uq}{{\underline{q}}}
\nc{\ur}{{\underline{r}}}
\nc{\us}{{\underline{s}}}\nc{\ut}{{\underline{t}}}
\nc{\uw}{{\underline{w}}}
\nc{\uz}{{\underline{z}}}
\nc{\ual}{{\underline{\alpha}}}\nc{\ualpha}{{\underline{\alpha}}}
\nc{\ugamma}{{\underline{\gamma}}}
\nc{\ula}{{\underline{\lambda}}}\nc{\umu}{{\underline{\mu}}}
\nc{\unu}{{\underline{\nu}}}\nc{\usigma}{{\underline{\sigma}}}
\nc{\utau}{{\underline{\tau}}}
\nc{\uN}{{\underline{N}}}\nc{\uM}{{\underline{M}}}
\nc{\uK}{{\underline{K}}}

\nc{\A}{{\mathfrak a}} \nc{\B}{{\mathfrak b}} \nc{\C}{{\mathfrak c}} 
\nc{\G}{{\mathfrak g}} \nc{\D}{{\mathfrak d}} \nc{\HH}{{\mathfrak h}}
\nc{\iii}{{\mathfrak i}}\nc{\mm}{{\mathfrak m}}\nc{\N}{{\mathfrak n}} 
\nc{\ttt}{{\mathfrak{t}}}\nc{\U}{{\mathfrak u}}\nc{\V}{{\mathfrak v}}

\nc{\SL}{{\mathfrak{sl}}}

\nc{\SG}{{\mathfrak S}}

\nc{\wt}{\widetilde} \nc{\wh}{\widehat}
\nc{\bn}{\begin{equation}}\nc{\en}{\end{equation}} \nc{\td}{\tilde}

%
%
%

\newcommand{\ldar}[1]{\begin{picture}(10,50)(-5,-25)
\put(0,25){\vector(0,-1){50}}
\put(5,0){\mbox{$#1$}} 
\end{picture}}

\newcommand{\luar}[1]{\begin{picture}(10,50)(-5,-25)
\put(0,-25){\vector(0,1){50}}
\put(5,0){\mbox{$#1$}}
\end{picture}}

\title[Universal algebras associated to Lie bialgebras]
{On some universal algebras associated to \\ the category of Lie 
bialgebras}

\author{B. Enriquez} 

\address{D\'epartement de Math\'ematiques et Applications, 
Ecole Normale Sup\'erieure, UMR 8553 du CNRS, 
45 rue d'Ulm, 75005 Paris, France}

\date{January 2001}

\begin{abstract} 
In our previous work (math/0008128), we studied the 
set $\Quant(\KK)$ of all universal quantization 
functors of Lie bialgebras over a field $\KK$ of 
characteristic zero, compatible with duals and doubles. 
We showed that $\Quant(\KK)$ is canonically isomorphic
to a product $\cG_0(\KK) \times \Sha(\KK)$, where $\cG_0(\KK)$
is a universal group and $\Sha(\KK)$ is a quotient set of a set 
$\cB(\KK)$ of families of Lie polynomials by the action of a group 
$\cG(\KK)$. We prove here that 
$\cG_0(\KK)$ is equal to the multiplicative group 
$1 + \hbar \KK[[\hbar]]$. So $\Quant(\KK)$ is `as close as 
it can be' to $\Sha(\KK)$. We also prove that the only
universal derivations of Lie bialgebras are multiples of the 
composition of the bracket with the cobracket.  Finally, we prove
that the stabilizer of any element of $\cB(\KK)$ is reduced to  
the $1$-parameter subgroup generated by the corresponding 
`square of the antipode'.    
\end{abstract}

\maketitle

\section{Main results}

\subsection{Results on $\Quant(\KK)$}

Let $\KK$ be a field of characteristic zero. 
In \cite{Enr}, we introduced the group $\cG_0(\KK)$
of all universal automorphisms of the adjoint representations of 
$\KK[[\hbar]]$-Lie bialgebras. 

Let us recall the definition of $\cG_0(\KK)$ 
more explicitly. Let $\hbar$ be a formal variable and 
let $\LBA_\hbar$ be the category of Lie bialgebras
over $\KK[[\hbar]]$, which are topologically free
$\KK[[\hbar]]$-modules. An element of $\cG_0(\KK)$
is a functorial assignment $(\A,[,],\delta_\A) \mapsto 
\rho_\A$, where for each object $(\A,[,],\delta_\A)$  of 
$\LBA_\hbar$, $\rho_\A$ is an element of $\End_{\KK[[\hbar]]}(\A)$, 
such that $(\rho_\A$ mod $\hbar) = \id_\A$, 
$\rho_{\A^*} = (\rho_\A)^t$, and for any $x,y$ in $\A$, 
$[\rho_\A(x),y] = \rho_\A([x,y])$, and $\rho_\A$
is given by a  composition of tensor products of the 
bracket and cobracket of $\A$, this composition being the 
same for all Lie bialgebras (we express the latter condition 
by saying that $\A\mapsto \rho_\A$ is universal).

View $1+\hbar\KK[[\hbar]]$ as a multiplicative subgroup 
of $\KK[[\hbar]]^\times$. 
There is a unique map $\al : 1 + \hbar \KK[[\hbar]]\to\cG_0(\KK)$,  
such that for any Lie bialgebra $\A$,
$(\al(\la))_\A = \la\id_\A$.  This map makes
$1+ \hbar \KK[[\hbar]]$ a subgroup of $\cG_0(\KK)$. 

We will show 

\begin{thm} \label{thm:G0}
$\cG_0(\KK)$ is equal to its subgroup $1 + \hbar \KK[[\hbar]]$. 
\end{thm}

In \cite{Enr}, we defined $\Quant(\KK)$ as the set of
all isomorphism classes of universal quantization functors of Lie bialgebras, 
compatible with duals and doubles. If we denote by $\LBA$ and
$\QUE$ the categories of Lie bialgebras and of quantized universal 
enveloping algebras over $\KK$, and by $\class : \QUE\to\LBA$
the semiclassical limit functor, then a 
universal quantization functor of Lie bialgebras, 
compatible with duals and doubles, 
is a functor $Q : \LBA\to \QUE$, such that 

1) $\class\circ Q$ is isomorphic to the identity; 

2) (universality) there exists an isomorphism 
of functors between $\A\mapsto Q(\A)$ and $\A\mapsto U(\A)[[\hbar]]$ 
(these are viewed as functors from $\LBA$ to the category of 
$\KK[[\hbar]]$-modules) with the following properties: 
if we compose this isomorphism with the symmetrisation map
$U(\A)[[\hbar]] \to S(\A)[[\hbar]]$, and if we transport the
operations of $Q(\A)$ on $S(\A)[[\hbar]]$, then the expansion in 
$\hbar$ of these operations yields maps 
$S^i(\A)\otimes S^j(\A)\to S^k(\A)$ and 
$S^i(\A) \to S^j(\A)\otimes S^k(\A)$; we require that these
maps be compositions of tensor products of the bracket and 
cobracket of $\A$, these compositions being independent of $\A$ 
(see \cite{Dr:open})

3) if $Q^\vee$ (resp., $D(Q)$) denotes the QUE-dual (resp., Drinfeld
double) of an object $Q$ of $\QUE$, and $D(\A)$ denotes the double
Lie bialgebra of an object $\A$ of $\LBA$, then there are canonical 
isomorphisms $Q(\A^*) \to Q(\A)^\vee$ and  
$Q(D(\A)) \to D(Q(\A))$; moreover, the universal $R$-matrix 
of $D(Q(\A))$ should be functorial (see \cite{Enr}). 

\medskip 

In \cite{Enr}, we also introduced
an explicit set $\Sha(\KK)$ of equivalence classes of 
families of Lie polynomials, satisfying associativity 
relations, and we constructed a canonical injection 
of $\Sha(\KK)$ in $\Quant(\KK)$. Moreover, we
constructed an action of $\cG_0(\KK)$ on $\Quant(\KK)$ and
showed that the map 
$$
\cG_0(\KK)\times \Sha(\KK) \to 
\Quant(\KK)
$$ 
given by the composition  
$\cG_0(\KK)\times \Sha(\KK) \subset \cG_0(\KK)\times
\Quant(\KK)\to \Quant(\KK)$ (in which the second map is the action map
of $\cG_0(\KK)$ on $\Quant(\KK)$) is a bijection. 
Theorem \ref{thm:G0} therefore implies

\begin{cor}
If $\A = (\A,[,],\delta_\A)$ is an object of $\LBA$ and 
$\la\in (1+\hbar \KK[[\hbar]])$, 
let $\A_\la$ be the object of $\LBA_\hbar$ isomorphic to 
$(\A,[,],\la\delta)$. The group $1+\hbar \KK[[\hbar]]$ acts 
freely on $\Quant(\KK)$ by the rule 
$(\la,Q)\mapsto Q_\la$, where $Q_\la$ is the functor
$\A\mapsto \wh Q(\A_\la)$ and $\wh Q$ is the natural extension of 
$Q$ to a functor from $\LBA_\hbar$ to $\QUE$. 
Then the map 
$$
\big( 1 + \hbar \KK[[\hbar]] \big) \times \Sha(\KK) \to \Quant(\KK) 
$$
given by the composition   
$\big( 1 + \hbar \KK[[\hbar]] \big)
\times \Sha(\KK) \subset \big( 1 + \hbar \KK[[\hbar]] \big) \times
\Quant(\KK)\to \Quant(\KK)$ is a bijection. 
\end{cor}

Therefore $\Quant(\KK)$ is `as close as it can be' to 
$\Sha(\KK)$. 

\subsection{Universal (co)derivations of Lie bialgebras}

Let $\cD$ (resp., $\cC$) be the space of all universal derivations 
(resp., coderivations) of Lie bialgebras. 
More explicitly, $\cD$ (resp., $\cC$) is the linear space of all functorial assignments 
$\A\mapsto\la_\A$, where for each object $\A$ of $\LBA$, $\la_\A$
belongs to $\End(\A)$, is universal in the above sense, and 
is a derivation (resp., coderivation) of the Lie algebra structure of $\A$. 
It is well-known that if $[,]_\A$ and $\delta_\A$ are the bracket and 
cobracket  maps of $\A$, then $[,]_\A \circ \delta_\A$ is a derivation of 
$\A$; e.g., if $\A$ is finite-dimensional, if $\sum_{i\in I} a_i \otimes b_i$
is the canonical element of $\A\otimes\A^*$ and if we set 
$u = \sum_{i\in I} [a_i,b_i]$,  then we have the identity
$([,]_\A \circ \delta_\A)(x) = [u,x]$ in the double Lie algebra 
of $\A$; and since $[,]_{\A^*}\circ\delta_{\A^*}$ is a derivation of 
$\A^*$, its transpose $[,]_\A\circ\delta_\A$ is a coderivation of 
$\A$. Then 

\begin{thm} \label{thm:ders}
$\cD$ and $\cC$ both coincide with
the one-dimensional vector space spanned by the 
assignment $\A\mapsto [,]_\A\circ\delta_\A$. 
\end{thm}

If $V$ is a vector space, we denote by 
$F(V)$ the free Lie algebra generated by $V$. Then the assignment
$\C\mapsto F(\C)$ is a functor from the category $\LCA$ of Lie coalgebras 
to $\LBA$. The proof of Theorem \ref{thm:ders} implies the following 
analogous statement
for the subcategory of $\LBA$ of free Lie 
algebras of Lie coalgebras. 

\begin{prop} \label{prop:ders:free}
Let $\C\mapsto\la_\C$ be a functorial assignment, where for each object 
$\C$ of $\LCA$, $\la_\C$ is a both a derivation and a 
coderivation of $F(\C)$. Then there exists a 
scalar $\la$, such that for any object $\C$ of $\LBA$, $\la_\C = 
\la [,]_{F(\C)}\circ\delta_{F(\C)}$. 
\end{prop}

\subsection{Isotropy of the action of $\cG(\KK)$ on $\cB(\KK)$}

We record here the definition of $\Sha(\KK)$. Let $\cB(\KK)$
be the set of families $(B_{pq})_{p,q\geq 0}$, such that 
for each $p,q$, $B_{pq}$ belongs to $FL_{p+q}[[\hbar]]$, 
$B_{10}(x) = B_{01}(x) = x$, $B_{p0} = B_{0p}$ if $p\neq 1$, 
$B_{11}(x,y) = [x,y]$, and for any integers $p,q,r$, the identity 
\begin{align*}
& \sum_{\al>0} 
\sum_{(p_\beta)_{\beta = 1,\ldots,\al}\in\Part_{\al}(p),
(q_\beta)_{\beta = 1,\ldots,\al}\in\Part_{\al}(q)} 
B_{\al r} \big( B_{p_1q_1}(x_1,\ldots,x_{p_1}|y_1,\ldots,y_{q_1}) 
\cdots \\ & \nonumber 
\cdots B_{p_\al q_\al}(
x_{\sum_{\beta = 1}^{\al - 1}p_\beta + 1},
\ldots, x_p|y_{\sum_{\beta = 1}^{\al - 1}q_\beta + 1},
\ldots, y_q)
|z_1,\ldots,z_r \big)
\\ & \nonumber 
= \sum_{\al>0} 
\sum_{(q_\beta)_{\beta = 1,\ldots,\al}\in\Part_{\al}(q),
(r_\beta)_{\beta = 1,\ldots,\al}\in\Part_{\al}(r)} 
B_{p \al} \big( x_1,\ldots, x_p| 
\\ & \nonumber 
B_{q_1r_1}(y_1,\ldots,y_{q_1}|z_1,\ldots,z_{r_1}) 
\cdots 
B_{q_\al r_\al}(
y_{\sum_{\beta = 1}^{\al - 1}q_\beta + 1},
\ldots, y_q|z_{\sum_{\beta = 1}^{\al - 1}r_\beta + 1},
\ldots, z_r)
\big)
\end{align*}
holds; here $FL_n$ is the multilinear part of the free Lie algebra
over $\KK$ with $n$ generators, and 
$\Part_\al(n)$ is the set of $\al$-partitions of $n$, 
i.e.\  the set of families $(n_1,\ldots,n_\al)$ of positive 
integers such that 
$n_1 + \cdots + n_\al = n$. Define $\cG(\KK)$ as the subset of $\prod_{n\geq 1}
FL_n[[\hbar]]$ of families $(P_n)_{n\geq 1}$, such that 
$P_1(x) = x$. Then we are going to define a group structure
on $\cG(\KK)$, and an action of $\cG(\KK)$ on 
$\cB(\KK)$; $\Sha(\KK)$ is the quotient set $\cB(\KK) / \cG(\KK)$.

Recall first that if 
$(\C,\delta_\C)$ is a Lie coalgebra and $B\in\cB(\KK)$, 
then there is a 
unique Hopf algebra structure on the completed tensor algebra 
$T(\C)[[\hbar]]$ with coproduct $\Delta^{\C}_{B} : 
T(\C)[[\hbar]] \to T(\C)^{\otimes 2}[[\hbar]]$ defined by 
$$
\Delta^{\C}_{B}(x) = x\otimes 1 + 1\otimes x
+ \sum_{p,q|p+q\geq 1} \hbar^{p+q - 1} \al_{pq}(\delta^{(B_{pq})}_\C(x))
$$
for any $x\in \C$, where for any $P$ in $FL_n[[\hbar]]$,
$\delta^{(P)}_\C$ is the map from 
$\C$ to $\C^{\otimes n}$ dual to the map from 
$(\C^*)^{\otimes n}$ to $\C^*$ defined by $P$
(when $\C$ is finite-dimensional), and $\al_{pq}$
is the map from $\C^{\otimes p+q}[[\hbar]]$ to 
$T(\C)\otimes T(\C)[[\hbar]]$ sending $x_1\otimes \cdots\otimes 
x_{p+q}$ to $(x_1\otimes\cdots\otimes x_p)\otimes
(x_{p+1}\otimes\cdots\otimes x_{p+q})$. 

If $P = (P_n)_{n\geq 1}$ belongs to $\cG(\KK)$, and 
$(\C,\delta_\C)$ is a Lie coalgebra,  
define $i^\C_P$ as the unique automorphism of 
$T(\C)[[\hbar]]$, such that for any $x$ in $\C$, we have
$$
i_P^\C(x) = x + \sum_{n\geq 2} \hbar^{n-1} \delta_\C^{(P_n)}(x). 
$$
Then the product $* : \cG(\KK)\times\cG(\KK) \to\cG(\KK)$ and
the operation $* : \cG(\KK)\times\cB(\KK) \to\cB(\KK)$ are uniquely
determined by the conditions that for any Lie coalgebra $(\C,\delta_\C)$, 
and any $P,Q$ in $\cG(\KK)$ and any $B$ in $\cB(\KK)$, we have 
$$
i^\C_{P*Q} = i^\C_P\circ i^\C_Q \quad\on{and}\quad  
\Delta^\C_{P*B} = 
(i^\C_P\otimes i^\C_P)\circ\Delta^\C_B\circ (i^\C_P)^{-1}. 
$$
Then one checks that for any $B$ in $\cB(\KK)$, there exists a unique
family $(S_n)_{n\geq 2}$, where $S_n\in FL_n[[\hbar]]$, such that 
for any Lie coalgebra $(\C,\delta_\C)$, the antipode $S_B^\C$ of 
the bialgebra $(T(\C)[[\hbar]],m_0,\Delta^\C_B)$ is such that 
for any $x\in\C$, 
$$
S^\C_B(x) = - x + \sum_{n|n\geq 2} \hbar^{n-1} \delta^{(S_n)}(x)
$$ 
($m_0$ is the multiplication map in $T(\C)[[\hbar]]$). 
It follows that there is a unique family $(\wt S_n)_{n\geq 2}$, 
where $\wt S_n\in FL_n[[\hbar]]$, such that 
$$
(S^\C_B)^2(x) = x + \sum_{n|n\geq 2} \hbar^{n-1} \delta^{(\wt S_n)}(x)
$$ 
for any $x\in\C$. Let us set $\wt S_1(x) = x$ and set 
$S_B^2 = (\wt S_n)_{n\geq 1}$. Then $S_B^2$ is an element of 
$\cG(\KK)$. Moreover, $\cG(\KK)$ is a pro-unipotent group. 
It follows that one can define the logarithm $\log(S_B^2)$
of $S_B^2$, and the corresponding one-parameter subgroup 
$\exp(\KK[[\hbar]]\log(S_B^2))$.  

Since for any Lie coalgebra $(\C,\delta_\C)$, we have 
$((S_B^C)^2\otimes (S_B^\C)^2)\circ \Delta_B^\C = 
\Delta_B^\C\circ (S_B^\C)^2$, we also have $S_B^2 * B = B$. 
It follows that for any element $g$ of 
$\exp(\KK[[\hbar]]\log(S_B^2))$, we have $g* B = B$,   
so $\exp(\KK[[\hbar]]\log(S_B^2))$ is contained in 
the isotropy group of $B$. 

\begin{prop} \label{isotropy}
For any $B$ in $\cB(\KK)$, the isotropy group of $B$ for 
the action of $\cG(\KK)$ on $\cB(\KK)$ is equal to  
$\exp(\KK[[\hbar]]\log(S_B^2))$. 
\end{prop}

\section{Proof of Theorem \ref{thm:G0}}

Let us define $\cE$ as the set of all universal 
$\KK[[\hbar]]$-module endomorphisms of Lie bialgebras. 
More explicitly, an element $\cE$ is a functorial 
assignment $(\A,[,],\delta_\A)\mapsto \eps_\A\in\End_{\KK[[\hbar]]}(\A)$, 
where $\A$ is an object of $\LBA_\hbar$ and 
the universality requirement means that 
$\eps_\A$ is given by a  composition of tensor products of the 
bracket and cobracket of $\A$, this composition being the 
same for all Lie bialgebras.  Then $\cE$ is a $\KK[[\hbar]]$-module, 
and $\cG_0(\KK)$ is a subset of $\cE$. 
We will first give a description of $\cE$ is terms of 
multilinear parts of free Lie algebras 
(Proposition \ref{prop:desc:E}). Then a computation in 
free algebras will prove Theorem \ref{thm:G0}. 

\subsection{Description of $\cE$} 
  
Define $FL_n$ as the multilinear part in each generator 
of the free Lie algebra over $\KK$ with $n$ generators.
Let $\SG_n$ act diagonally on $FL_n\otimes FL_n$ by 
simultaneous permutation of the generators $x_1,\ldots,x_n$
and $y_1,\ldots,y_n$ of each factor.    

We define a linear map $p\mapsto (\A\mapsto i(p)_\A)$ from  
$\wh\oplus_{n|n\geq 1}  (FL_n\otimes
FL_n)_{\SG_n}[[\hbar]]$ to $\cE$ as follows.  

If $Q$ is an element of $FL_n$, which we write
$Q = \sum_{\sigma\in\SG_n} Q_\sigma 
x_{\sigma(1)} \cdots x_{\sigma(n)}$, then we set 
$$
\delta_\A^{(Q)}(x) = {1\over n} \sum_{\sigma\in\SG_n}
Q_\sigma \big(
(\id_\A^{\otimes n-1}\otimes \delta_\A)\circ 
\cdots\circ\delta_\A(x)\big)^{(\sigma(1) \ldots \sigma(n))}
$$ 
for any $x\in\A$. 

If $p = \sum_\al P_\al\otimes Q_\al$ is an element of
$FL_n\otimes FL_n$, define $i(p)_\A$ as the endomorphism of 
$\A$ such that 
\begin{equation} \label{def:map}
i(p)_\A(x) = \sum_\al P_\al(\delta^{(Q_\al)}(x))
\end{equation}
for any $x\in\A$. 
This  maps factors through a linear map $p\mapsto i(p)_\A$
from $(FL_n\otimes FL_n)_{\SG_n}$ to $\End_{\KK[[\hbar]]}(\A)$, 
and induces a linear map $p\mapsto \wh i(p)_\A$  from 
$\wh\oplus_{n|n\geq 1} (FL_n\otimes FL_n)_{\SG_n}[[\hbar]]$ 
to $\End_{\KK[[\hbar]]}(\A)$ (here $\wh\oplus$ is the $\hbar$-adically 
completed direct sum).  

Then if $\A$ is finite-dimensional over $\KK[[\hbar]]$, 
and if we express the canonical element of $\A\otimes\A^*$
as $\sum_{i\in I} a_i\otimes b_i$, then 
$$
i(p)_\A(x) = \sum_\al\sum_{i_1,\ldots,i_n\in I}
\langle x, Q_\al(b_{i_1},\ldots,b_{i_n}) \rangle
P_\al(a_{i_1},\ldots,a_{i_n}) . 
$$
 
\begin{prop} \label{prop:desc:E}
The linear map $\wh i$ from 
$\wh\oplus_{n|n\geq 1}  (FL_n\otimes
FL_n)_{\SG_n}[[\hbar]]$ to $\cE$ defined by $p\mapsto 
(\A\mapsto \wh i(p)_\A)$ 
is a linear isomorphism.
\end{prop}

{\em Proof.} If $n$ and $m$ are $\geq 1$, define 
$\cE_{n,m}$ as the vector space of
all universal linear homomorphisms from $\A^{\otimes n}$ to 
$\A^{\otimes m}$ (`universal' again means that these homomorphisms
are compositions of tensor products of the bracket and cobracket 
map, this composition being the same for each $\A$). 
Then $\cE$ is just $\cE_{1,1}$. 
The direct sum $\wh\oplus_{n,m|n,m\geq 1}\cE_{n,m}$ may be defined 
formally as the smallest $\hbar$-adically complete
vector subspace of the space of all 
functorial assignments $\A\mapsto \rho_\A\in
\wh\oplus_{n,m|n,m\geq 1}\Hom_{\KK[[\hbar]]}(\A^{\otimes n},
\A^{\otimes m})$, containing the assignments 
$\A\mapsto\id_\A\in\Hom_{\KK[[\hbar]]}(\A,\A)$, 
the bracket and the cobracket 
operations, stable under the external 
tensor products operations 
$\Hom_{\KK[[\hbar]]}(\A^{\otimes n},\A^{\otimes m}) 
\otimes 
\Hom_{\KK[[\hbar]]}(\A^{\otimes n'},\A^{\otimes m'})
\to \Hom_{\KK[[\hbar]]}(\A^{\otimes n+n'},\A^{\otimes m+m'})$, 
under the natural actions of the symmetric groups $\SG_n$ 
and $\SG_m$ on $\Hom_{\KK[[\hbar]]}(\A^{\otimes n},\A^{\otimes m})$, 
and under the composition operation. 

Define $\cE_{n,m}^{finite}$ and $\cE^{finite}$ as the analogues
of $\cE_{n,m}$ and $\cE$, where the category of Lie 
bialgebras is replaced by the category of finite-dimensional
(over $\KK[[\hbar]]$) $\KK[[\hbar]]$-Lie bialgebras.  
Then restriction to this subcategory of $\LBA_\hbar$
induces linear maps $\cE_{n,m}\to\cE_{n,m}^{finite}$
and  $\wh\oplus_{n,m|n,m\geq 1}\cE_{n,m}\to\wh\oplus_{n,m|n,m\geq 1}
\cE_{n,m}^{finite}$. 

Let us define $F^{(n,m)}$ as follows 
\begin{align*}
& 
F^{(n,m)} = 
\\ & 
\bigoplus_{(p_{ij})\in\NN^{\{1,\ldots,n\}\times\{1,\ldots,p\}}}
\big( 
\bigotimes_{i = 1}^n FL_{\sum_{j = 1}^m p_{ij}} 
\otimes
\bigotimes_{j = 1}^m FL_{\sum_{i = 1}^n p_{ij}}
\big)_{\prod_{(i,j)\in\{1,\ldots,n\}\times\{1,\ldots,m\}} \SG_{p_{ij}}} ; 
\end{align*}
the generators of the $i$th factor of the first tensor product are
$x^{(ij)}_\al$, $j = 1,\ldots,m$, $\al = 1,\ldots,p_{ij}$,  
and the generators of the $j$th factor of the second tensor
product are $y^{(ij)}_\al$, $i = 1,\ldots,n$, $\al = 1,\ldots,p_{ij}$; 
the group $\SG_{p_{ij}}$ acts by simultaneously permuting the 
generators $x^{(ij)}_\al$ and $y^{(ij)}_\al$, $\al\in\{1,\ldots,p_{ij}\}$. 

Then there is a unique linear map $i_{n,m} : 
F^{(n,m)}\to\cE^{finite}_{m,n}$, 
such that if 
\begin{align*}
& p = 
\\ & \sum_\la \bigotimes_{i = 1}^n P^\la_i(x^{(ij)}_\al; 
j = 1,\ldots,m; \al = 1,\ldots,p_{ij})
\otimes \bigotimes_{j = 1}^m Q^\la_j(y^{(ij)}_\al; 
i = 1,\ldots,n; \al = 1,\ldots,p_{ij}) ,
\end{align*}
and $x_1,\ldots,x_m$ belong to a Lie bialgebra $\A$, then 
\begin{align} \label{def:inm}
& (i_{n,m}(p))_\A(x_1\otimes\cdots\otimes x_m) = 
\\ & \nonumber
\sum_\la \sum_{i^{(11)}_1\in I,\ldots,i^{(nm)}_{p_{nm}}\in I}
\prod_{j = 1}^m \langle x_j, Q^\la_j(b(i^{(ij)}_\al); 
i = 1,\ldots,n: \al = 1,\ldots,p_{ij})\rangle  
\\ & \nonumber
\bigotimes_{i = 1}^n P^\la_i(a(i^{(ij)}_\al); j = 1,\ldots,m;  
\al = 1,\ldots,p_{ij}) . 
\end{align}
Here we write the canonical element of $\A\otimes\A^*$ in 
the form $\sum_{i\in I} a(i)\otimes b(i)$. 
Formula (\ref{def:inm}) has an obvious generalization 
when $\A$ is an arbitrary Lie bialgebra. This means
that $i_{n,m}$ factors through a map (also denoted $i_{n,m}$) 
from $F^{(n,m)}$ to $\cE_{m,n}$; this map induces linear
maps $\wh i_{n,m} : F^{(n,m)}[[\hbar]]\to \cE_{n,m}$ and 
$\wh\oplus_{n,m|n,m\geq 1} \wh i_{n,m} : 
\wh\oplus_{n,m|n,m\geq 1} F^{(n,m)}[[\hbar]]\to 
\wh\oplus_{n,m|n,m\geq 1}\cE_{n,m}$.   

Let us show that $\wh\oplus_{n,m|n,m\geq 1} \wh i_{n,m}$ is 
surjective. 
For this, let us study $\wh\oplus_{n,m\geq 1} \Imm(\wh i_{n,m})$.  
This is a subspace of the space of all functorial 
assignments 
$$
\A\mapsto \rho_\A\in\wh\oplus_{n,m|n,m\geq 1}
\Hom_{\KK[[\hbar]]}(\A^{\otimes n},\A^{\otimes m}).
$$ 
Let us show that it shares all the properties  of 
$\oplus_{n,m|n,m\geq 0}\cE_{n,m}$. The identity is the image of 
the element $x^{(11)}_1\otimes y^{(11)}_1$ of $F^{(1,1)}$, 
the bracket is the 
image of the element $[x^{(11)}_1,x^{(12)}_1]\otimes 
y^{(11)}_1 \otimes y^{(12)}_1$ of $F^{(1,2)}$, and the 
cobracket is the image of the element  
$x^{(11)}_1 \otimes x^{(21)}_1\otimes 
[y^{(11)}_1,y^{(21)}_1]$ of $F^{(2,1)}$. The fact that  
$\wh\oplus_{n,m\geq 0} \Imm(\wh i_{n,m})$ is stable under the 
composition follows from the following Lemma. 

\begin{lemma}
Let $P$ and $Q$ be Lie polynomials in $FL_n$ and $FL_m$
respectively. Then there exist an element 
$p = \sum_{\up\in \NN^{ \{1,\ldots,n\}\times\{1,\ldots,m\} } }
\sum_{\la} (\otimes_{i = 1}^n P^{\up,\la}_i)\otimes 
(\otimes_{j = 1}^m Q^{\up,\la}_j )$ of $F^{(n,m)}$, 
such that if $\A$ is any 
finite-dimensional Lie bialgebra over $\KK$, and 
$\sum_{i\in I} a(i)\otimes b(i)$ 
is the canonical element of $\A\otimes \A^*$, then 
\begin{align*}
& 
\sum_{i_1,\ldots,j_m\in I} 
\langle Q(a(j_1),\ldots,a(j_m)) , P(b(i_1),\ldots,b(i_n))\rangle
(\otimes_{\al = 1}^n a(i_\al)) \otimes
(\otimes_{\beta = 1}^m b(j_\beta)) 
\\ & = 
\sum_{\up\in \NN^{ \{1,\ldots,n\}\times\{1,\ldots,m\} } }
\sum_{\la} \sum_{i^{(11)}_1\in I,\ldots,i^{(nm)}_{p_{nm}}\in I}
(\bigotimes_{i = 1}^n P^{\up,\la}_i( a(i^{(ij)}_\al); 
j = 1,\ldots,m; \al = 1,\ldots,p_{ij} ) )
\\ & 
\otimes 
(\bigotimes_{j = 1}^m Q^{\up,\la}_j (b(i^{(ij)}_\al); 
i = 1,\ldots,n; \al = 1,\ldots,p_{ij}) ) . 
\end{align*}
\end{lemma}

{\em Proof of Lemma.}
We may assume that $P$ and $Q$ have the form 
$P(x_1,\ldots,x_n) = [x_1,[x_2,\ldots,x_n]]$ and
$Q(y_1,\ldots,y_m) = [y_1,[y_2,\ldots,y_m]]$.
Then the invariance of the canonical bilinear form in 
$D(\A)$ and the fact that $\sum_{i\in I} a(i)\otimes b(i)$
satisfies the classical Yang-Baxter identity in $D(\A)$ 
imply the following formula. 
If $\al$ is an integer and $\uk = (k_1,\ldots,k_\al)$ is 
a sequence of integers such that $1\leq k_1 < \ldots < k_\al < m$, 
let $\uk' = (k'_1,\ldots,k'_m)$ be the sequence such that 
$k'_i = k_i$ for $i = 1,\ldots,\al$, $(k'_i)_{i>\al}$ is decreasing
and $\{k'_1,\ldots,k'_m\} = \{1,\ldots,m\}$. 
Then  
\begin{align*}
& \sum_{i_1,\ldots,j_m\in I} 
\langle [a(j_1),[a(j_2),\ldots,a(j_m)]],  
[b(i_1),[b(i_2),\ldots,b(i_n)]] \rangle 
\\ & 
a(i_1)\otimes\cdots\otimes 
a(i_n) \otimes b(j_1)\otimes\cdots\otimes b(j_m) 
\\ &  
= \sum_{\al = 1}^{m-1} \sum_{k_1,\ldots,k_\al | 1\leq k_1 < \cdots 
< k_\al < m} \sum_{s = 0}^{m - \al} \kappa(s) \sum_{i_1,\ldots,j_m\in I}
\\ & 
\langle 
[a(j_{k'_m}),\ldots,[a(j_{k'_{\al+s+2}}),a(j_{k'_{\al+s+1}})]] , 
[b(i_2),[b(i_3),\ldots,b(i_n)]]\rangle 
\\ & 
[a(j_{k'_1}),\cdots [a(j_{k'_{\al+s}}),a(i_1)]]
\otimes a(i_2)\otimes\cdots\otimes a(i_m)
\\ & 
\otimes b(j_1)\otimes\cdots\otimes [b(j_{k'_{\al +s} + 1}),b(i_1)]
\otimes\cdots\otimes b(j_m) ,  
\end{align*}  
where $\kappa(0) = -1$ and $\kappa(s) = (-1)^s$ if $s\neq 0$.  
The Lemma then follows by induction on $n$ and $m$. 
\hfill \qed\medskip

{\em End of proof of Proposition.}
The other properties of  $\wh\oplus_{n,m|n,m\geq 0}\cE_{n,m}$ 
are obviously shared by 
$\wh\oplus_{n,m\geq 0} \Imm(\wh i_{n,m})$. 
So $\wh\oplus_{n,m\geq 0} \Imm(\wh i_{n,m})$ is contained in 
$\wh\oplus_{n,m|n,m\geq 0}\cE_{n,m}$ and shares all its properties; 
since $\wh\oplus_{n,m|n,m\geq 0}\cE_{n,m}$ is the smallest vector
subspace of the space of functorial assignments 
$\A\mapsto \rho_\A\in\wh\oplus_{n,m|n,m\geq 1} \Hom_{\KK[[\hbar]]}
(\A^{\otimes n},\A^{\otimes m})$ with these
properties, we obtain  
$\wh\oplus_{n,m|n,m\geq 1} \Imm(\wh i_{n,m}) 
= \wh\oplus_{n,m|n,m\geq 1}\cE_{n,m}$. 
This proves that $\wh\oplus_{n,m|n,m\geq 1} \wh i_{n,m}$ is 
surjective. Since $\wh i = \wh i_{1,1}$, this implies that 
$\wh i$ is surjective.  

\medskip 

Let us now show that the map $\wh i$
is injective. Let $(p_n)_{n\geq 1}$ be a family such that  
$p_n\in (FL_n\otimes FL_n)_{\SG_n}[[\hbar]]$ and 
$\wh i(\sum_{n|n\geq 1} p_n) = 0$. Then if $\A$ is any 
Lie bialgebra, then $\sum_{n|n\geq 1}\wh i(p_n)_\A = 0$. 

If $V$ is any vector space, let us denote by $F(V)$
the free Lie algebra generated by $V$. If $(\C,\delta_\C)$
is a Lie coalgebra, then the map $\C\to\wedge^2 F(\C)$
defined as the composition of $\C\to\wedge^2 \C\to \wedge^2 F(\C)$
of the cobracket map of $\C$ with the canonical inclusion
extends to a unique cocycle map $\delta_{F(\C)} : 
F(\C)\to\wedge^2 F(\C)$. 
Then $(F(\C),[,],\delta_{F(\C)})$ is a Lie bialgebra. 
The assignment $\C\mapsto F(\C)$ is a functor from the category 
$\LCA$ of Lie coalgebras to $\LBA$. Then if $(\C,\delta_\C)$
is any Lie coalgebra, we have
\begin{equation} \label{MA}
\sum_{n|n\geq 1} \wh i(p_n)_{F(\C)} = 0 . 
\end{equation}

Define $FA_n$ as the multilinear part of the free 
algebra with generators $x_1,\ldots,x_n$. 
Then $\SG_n$ acts on $FA_n\otimes FL_n$
by simultaneously permuting the generators $x_1,\ldots,x_n$
of $FA_n$ and $y_1,\ldots,y_n$ of $FL_n$. 
The injection $FL_n\subset FA_n$ induces a linear
map $(FL_n\otimes FL_n)_{\SG_n}\to (FA_n\otimes FL_n)_{\SG_n}$;  
since $\SG_n$ is finite, this linear
map is an injection. Moreover, the map $FL_n \to 
(FA_n\otimes FL_n)_{\SG_n}$, sending $P$ to the class of
$x_1\cdots x_n \otimes P(y_1,\ldots,y_n)$, is a 
linear isomorphism.

For each $n$, define $\bar p_n$ as the element of $FL_n$
such that the equality 
$$
p_n = x_1\cdots x_n\otimes \bar p_n(y_1,\ldots,y_n)
$$  
holds in $(FA_n\otimes FL_n)_{\SG_n}[[\hbar]]$. 

The restriction of $\wh i(p_n)_{F(\C)}$ to $\C\subset F(\C)$
is a linear map $\wh i(p_n)_{F(\C)|\C}$ 
from $\C$ to $F(\C)[[\hbar]]$. The image of this map 
is actually contained in the degree $n$ part $F(\C)_n[[\hbar]]$ of
$F(\C)[[\hbar]]$. The space $F(\C)_n$ is a vector subspace of 
$\C^{\otimes n}$. 
Moreover, formula (\ref{def:map}) shows that the composition of  
$\wh i(p_n)_{F(\C)|\C}$ with the canonical inclusion
$F(\C)_n[[\hbar]]\subset \C^{\otimes n}[[\hbar]]$ coincides with $\delta^{(\bar p_n)}$.  
So if $(\C,\delta_\C)$ is any Lie coalgebra, the map 
$\sum_{n|n\geq 0}\delta^{(\bar p_n)} : \C \to \oplus_{n|n\geq 0}
\C^{\otimes n}[[\hbar]]$ is zero. 
Now the linear map $FL_n \to \{$functorial assignments 
$\C\mapsto \tau_\C\in\Hom(\C,\C^{\otimes n})$, where $\C$ is an 
object of $\LCA\}$ defined by $P\mapsto \delta^{(P)}$, 
is injective. This 
implies that each $\bar p_n$ is zero. So $\wh i$ is 
injective. 

This ends the proof of Proposition \ref{prop:desc:E}. 
\hfill \qed\medskip

\begin{remark} \label{rem:gon}
More generally, one may show that $\cE_{n,m}$ 
is isomorphic to $F^{(m,n)}$.  
\end{remark}

\begin{remark}
It would be interesting to 
understand 1) the algebra structure of $F^{(1,1)}$ 
provided by the isomorphism of Proposition \ref{prop:desc:E}
and 2) the algebra structure of $\oplus_{n,m|n,m\geq 1}
F^{(n,m)}$ provided by Remark \ref{rem:gon}.  As we noted before, 
the latter algebra is also equipped with natural 
operations of the symmetric groups $\SG_n$ and $\SG_m$
on each component $F^{(n,m)}$, and external product maps 
 $F^{(n,m)} \otimes F^{(n',m')} \to F^{(n+n',m+m')}$.  
\end{remark}

\subsection{Proof of Theorem \ref{thm:G0}}

Let $(\A\mapsto \rho_\A)$ belong to $\cG_0(\KK)$. 
Recall that this means that $(\A\mapsto \rho_\A)$
belongs to $\cE$, in particular, $\rho_\A$
belongs to $\End_{\KK[[\hbar]]}(\A)$ for any object 
$\A$ of $\LBA_\hbar$. The condition that $(\A\mapsto \rho_\A)$ 
belongs to $\cG_0(\KK)$ implies that $\rho_\A$ satisfies the identity
$\rho_\A([x,y]) = [\rho_\A(x),y]$ for any $x,y$ in $\A$. 
Let $p$ be the preimage of $\rho_\A$ by the map $\wh i$. 
Then there is a unique sequence $(p_n)_{n\geq 1}$, where 
$p_n$ belongs to $(FL_n\otimes FL_n)_{\SG_n}[[\hbar]]$, 
such that $p = \sum_{n|n\geq 1} p_n$. Then if we set 
$p_n = \sum_\al P^{(n)}_\al\otimes Q^{(n)}_\al$, and if 
$\A$ is finite-dimensional, and $\sum_{i\in I} a(i)\otimes b(i)$
is the canonical element of $\A\otimes\A^*$, then we have 
$$
\rho_\A(x) = \sum_{n|n\geq 1} \sum_{i_1,\ldots,i_n\in I}
\langle x, Q^{(n)}_\al (b(i_1),\ldots, b(i_n)) \rangle 
P^{(n)}_\al (a(i_1),\ldots,a(i_n)) . 
$$
Then 
$$
\rho_\A([x,y]) = \sum_{n|n\geq 1} \sum_{i_1,\ldots,i_n\in I}
\langle [x,y], Q^{(n)}_\al (b(i_1),\ldots, b(i_n)) \rangle 
P^{(n)}_\al (a(i_1),\ldots,a(i_n)) . 
$$

\begin{lemma} \label{lemma:r}
If $\xi$ belongs to $\A^*$, $[\sum_{i\in I} a(i)\otimes b(i), 
\xi\otimes 1 + 1\otimes \xi]$ belongs to $\A^*\otimes\A_*$, and 
we have 
$$
\langle  [x,y], \xi \rangle =
\langle x\otimes y, [\sum_{i\in I} a(i)\otimes b(i), 
\xi\otimes 1 + 1\otimes \xi] \rangle. 
$$
\end{lemma}

{\em Proof of Lemma.}
The first statement follows from the fact that 
$\sum_i a(i)\otimes b(i) + b(i)\otimes a(i)$ is $D(\A)$-invariant.

Let us prove the second statement. 
The invariance of the bilinear form of $D(\A)$ implies
that  
$\langle  [x,y], \xi \rangle$ is equal to 
$\langle x, [y,\xi]\rangle$. This is equal to 
$\langle \sum_{i\in I} [a(i),\xi]\otimes b(i), x\otimes y\rangle$. 
Since $\A$ is an isotropic subspace of $D(\A)$, this is the
same as  
$$
\langle \sum_{i\in I} [a(i),\xi]\otimes b(i) 
+ a(i)\otimes [b(i),\xi], x\otimes y\rangle.
$$ 
\hfill \qed\medskip 

So we get 
\begin{align} \label{motti}
& \rho_\A([x,y]) = \sum_{n|n\geq 1} \sum_\al \sum_{i_1,\ldots,i_n\in I} 
\\ & \nonumber 
\langle x \otimes y, \sum_i 
[a(i), Q^{(n)}_\al (b(i_1),\ldots, b(i_n))]\otimes b(i) 
+ a(i) \otimes [b(i), Q^{(n)}_\al (b(i_1),\ldots, b(i_n))] \rangle 
\\ & \nonumber 
P^{(n)}_\al (a(i_1),\ldots,a(i_n)) . 
\end{align}
On the other hand, we have 
$$
[\rho_\A(x),y] = \sum_{n|n\geq 0} \sum_{i_1,\ldots,i_n\in I}
\langle x, Q_\al^{(n)}(b(i_1),\ldots,b(i_n))\rangle 
[P_\al^{(n)}(a(i_1),\ldots, a(i_n)),y]
$$ 
so 
\begin{align} \label{gur}
& [\rho_\A(x),y] = 
\\ & \nonumber 
\sum_{n|n\geq 0} \sum_\al  \sum_{i_1,\ldots,i_n\in I}
\sum_{i\in I}
\langle x\otimes y, Q_\al^{(n)}(b(i_1),\ldots,b(i_n)) \otimes b(i)\rangle 
[P_\al^{(n)}(a(i_1),\ldots, a(i_n)),a(i)] .  
\end{align}
Comparing (\ref{motti}) and (\ref{gur}), and using 
the first part of Lemma \ref{lemma:r}, we get 
\begin{align} \label{bar:lev}
& \sum_{n|n\geq 1} \sum_\al \sum_{i_1,\ldots,i_n,i\in I} 
[P^{(n)}_\al(a(i_1),\ldots,a(i_n)),a(i)]\otimes 
Q^{(n)}_\al(b(i_1),\ldots,b(i_n)) 
\otimes b(i)
\\ & \nonumber 
= \sum_{n|n\geq 1} \sum_\al \sum_{i_1,\ldots,i_n,i\in I} 
P^{(n)}_\al(a(i_1),\ldots,a(i_n)) \otimes a(i) \otimes 
[b(i),Q^{(n)}_\al (b(i_1),\ldots,b(i_n))]
\\ & \nonumber
+ P^{(n)}_\al(a(i_1),\ldots,a(i_n)) \otimes 
[a(i),Q^{(n)}_\al(b(i_1),\ldots,b(i_n))] \otimes b(i) . 
\end{align}
In fact, it is easy to see that the identity 
\begin{align} \label{univ:bar:lev}
& \sum_{n\geq 1} \sum_\al  
[P^{(n)}_\al(x_1,\ldots,x_n),x]\otimes 
Q^{(n)}_\al(y_1,\ldots,y_n) 
\otimes y
\\ & \nonumber 
= \sum_{n|n\geq 1} \sum_\al  
P^{(n)}_\al(x_1,\ldots,x_n) \otimes x \otimes 
[y,Q^{(n)}_\al (y_1,\ldots,y_n)]
\\ & \nonumber
+ P^{(n)}_\al(x_1,\ldots,x_n) \otimes 
[x,Q^{(n)}_\al(y_1,\ldots,y_n)] \otimes y  
\end{align}
holds in $F^{(abb)}$
(in the notation of \cite{Enr}), which is the `Lie part'
of a universal algebra for solutions of the classical 
Yang-Baxter equation (CYBE). Equation (\ref{bar:lev}) is then 
a consequence of (\ref{univ:bar:lev}).  

Applying the Lie bracket to the two last tensor 
factors of (\ref{bar:lev}), we obtain 
\begin{align} \label{tam} 
& \sum_{n|n\geq 1} \sum_\al \sum_{i_1,\ldots,i_n,i\in I}
[P^{(n)}_\al(a(i_1),\ldots,a(i_n)),a(i)]\otimes 
[Q^{(n)}_\al(b(i_1),\ldots,b(i_n)) , b(i)]
\\ & \nonumber 
= \sum_{n|n\geq 1}  \sum_\al \sum_{i_1,\ldots,i_n,i\in I} 
P^{(n)}_\al(a(i_1),\ldots,a(i_n)) \otimes [[a(i),b(i)],
Q^{(n)}_\al (b(i_1),\ldots,b(i_n))] . 
\end{align}
Since $\sum_{i\in I} a(i)\otimes b(i)$ satisfies CYBE, 
we have 
$$
\sum_{j\in I} a(j) \otimes [\sum_{i\in I} [a(i),b(i)],b(j)] = 
\sum_{i,j\in I} [a(i),a(j)]\otimes[b(i),b(j)]. 
$$ 
So identity (\ref{tam}) is rewritten as  
\begin{align} \label{senora} 
& \sum_{n|n\geq 1} \sum_\al \sum_{i_1,\ldots,i_n,i\in I}
[P^{(n)}_\al(a(i_1),\ldots,a(i_n)),a(i)]\otimes 
[Q^{(n)}_\al(b(i_1),\ldots,b(i_n)) , b(i)]
\\ & \nonumber 
= \sum_{n|n\geq 1}  \sum_\al \sum_{k = 1}^n 
\sum_{i_1,\ldots,i_n,i\in I}
P^{(n)}_\al(a(i_1),\ldots,[a(i),a(i_k)],\ldots,a(i_n)) 
\\ & \nonumber 
\otimes Q^{(n)}_\al (b(i_1),\ldots,[b(i),b(i_k)]\ldots,b(i_n)) . 
\end{align}
On the other hand, one easily derives from 
(\ref{univ:bar:lev}) the identity 
\begin{align} \label{nasi} 
& \sum_{n|n\geq 1} \sum_\al
[P^{(n)}_\al(x_1,\ldots,x_n),x]\otimes 
[Q^{(n)}_\al(y_1,\ldots,y_n),y]
\\ & \nonumber 
= \sum_{n|n\geq 1} \sum_{k = 1}^n \sum_\al 
P^{(n)}_\al(x_1,\ldots,[x,x_k],\ldots,x_n) \otimes 
Q^{(n)}_\al (y_1,\ldots,[y,y_k]\ldots,y_n)  
\end{align}
valid in $F^{(1,1)}$; this identity 
is the universal version of (\ref{senora}).   
Separating homogeneous components, we get for each 
$n\geq 1$
\begin{align} \label{mendez} 
& \sum_\al  
[P^{(n)}_\al(x_1,\ldots,x_n),x]\otimes 
[Q^{(n)}_\al(y_1,\ldots,y_n),y]
\\ & \nonumber 
= \sum_{\al} \sum_{k = 1}^n 
P^{(n)}_\al(x_1,\ldots,[x,x_k],\ldots,x_n) \otimes 
Q^{(n)}_\al (y_1,\ldots,[y,y_k]\ldots,y_n)  . 
\end{align}
For each $n\geq 1$, let $R_n$ be the element of $FL_n[[\hbar]]$ 
such that 
the identity $R_n(x_1,\ldots,x_n)\otimes y_1\cdots y_n
= \sum_\al P^{(n)}_\al(x_1,\ldots,x_n) \otimes Q^{(n)}_\al
(y_1,\ldots,y_n)$ holds in $(FL_n\otimes FA_n)_{\SG_n}[[\hbar]]$. 
Then (\ref{mendez}) implies that $R_n$ satisfies 
identity 
\begin{equation} \label{cordovero}
[R_n(x_1,\ldots,x_n),x_{n+1}] - [R_n(x_2,\ldots,x_{n+1}),x_1]
= 2 \sum_{k = 1}^n R_n(x_1,\ldots,[x_k,x_{k+1}],\ldots,x_{n+1})
\end{equation}
in $FL_{n+1}[[\hbar]]$. 

There are unique elements $R^{(i)}_n$ of $FA_{n-1}[[\hbar]]$
($i = 1,\ldots,n$), such that 
$$
R_n(x_1,\ldots,x_n) = \sum_{i = 1}^n x_i R_n^{(i)}(x_1,\ldots,
x_{i-1},x_{i+1},\ldots,x_n). 
$$
Let us view (\ref{cordovero}) as an identity in $FA_n[[\hbar]]$, 
and let us project it on $\oplus_{\sigma| \sigma\in \SG_{n}, 
\sigma(1) = 1}$ $\KK[[\hbar]] x_{\sigma(1)} \cdots x_{\sigma(n)}$
parallel to $\oplus_{\sigma| \sigma\in \SG_{n}, 
\sigma(1) \neq 1} \KK[[\hbar]] x_{\sigma(1)} \cdots x_{\sigma(n)}$. 
This means that we select in this identity the terms `starting with $x_1$'. 
This yields 
\begin{align} \label{luria}
& R_n^{(1)}(x_2,\ldots,x_n) x_{n+1} + R_n(x_2,\ldots,x_{n+1})
\\ & \nonumber 
= 2 x_2 R_n^{(1)}(x_3,\ldots,x_{n+1}) 
+ 2\sum_{i = 2}^n R_n^{(1)}(x_2,\ldots,[x_i,x_{i+1}],\ldots,x_{n+1}). 
\end{align}
This is an equality in $FA_n[[\hbar]]$. Let us denote by 
$\cF\cA_n$ the free $\KK[[\hbar]]$-algebra with generators 
$x_2,\ldots,x_{n+1}$, and by $\cF\cL_n$ the free Lie 
algebra with the same generators. Then 
$FA_n[[\hbar]] \subset \cF\cA_n$, and $\cF\cA_n$
is the universal enveloping algebra $U(\cF\cL_n)$
of $\cF\cL_n$. This structure of an enveloping 
algebra defines a filtration on $\cF\cA_n$.
An element of $\cF\cA_n$ has degree $\leq \beta$
for this filtration iff it can be expressed as a 
polynomial of degree $\leq\beta$ in elements of 
$\cF\cL_n$. 

Let $\al$ be the degree of $R_n^{(1)}$ for the 
analogous filtration of $\cF\cA_{n-1}$. Let us assume that 
$\al>0$. Then 
$R_n^{(1)}(x_2,\ldots,x_n) x_{n+1}$ and
$x_2 R_n^{(1)}(x_3,\ldots,x_{n+1})$
both have degree $\al +1$ in $\cF\cA_n$; 
$R_n^{(1)}(x_2,\ldots,[x_i,x_{i+1}],\ldots,x_{n+1})$ 
has degree $\al$ in $\cF\cA_n$; and 
$R_n(x_2,\ldots,x_{n+1})$ has degree $1$
(its degree is $\leq 1$, but if this degree is zero, then 
$R_n$ vanishes identically). 

Let $\bar R_n^{(1)}$ be the image of $R_n^{(1)}$
in the associated graded of $\cF\cA_n$, which 
is the symmetric algebra $S(FL_n)[[\hbar]]$ of $FL_n$.  
Then (\ref{luria}) implies the identity 
\begin{equation} \label{sholem}
\bar R_n^{(1)}(x_2,\ldots,x_n) x_{n+1} 
= 2 x_2 \bar R_n^{(1)}(x_3,\ldots,x_{n+1}) 
\end{equation}
in $S(\cF\cL_n)[[\hbar]]$. Recall that in this
identity, $\bar R_n^{(1)}(y_1,\ldots,y_{n-1})$ is a polynomial 
in variables $P_{k,\al}(y_{i_1},\ldots,y_{i_k})$, where
$k$ runs over $1,\ldots,n-1$, 
$i_1,\ldots,i_k$ runs over all sequences of integers
such that $1\leq i_1 < \cdots < i_k\leq n-1$, and 
$P_{k,\al}$ runs over a basis of $FL_k$, so it is a polynomial in 
$\sum_{k = 1}^{n-1} \pmatrix n-1 \\ k \endpmatrix (k-1)!$ variables.
Moreover, if $(\delta_1,\ldots,\delta_{n-1})$ is the canonical 
basis of $\NN^n$, and we say that the variables 
$P_{k,\al}(y_{i_1},\ldots,y_{i_k})$ have multidegree 
$\delta_{i_1} + \cdots + \delta_{i_k}$, then 
$\bar R_n^{(1)}$ is homogeneous of multidegree 
$\delta_1 + \cdots + \delta_{n-1}$.  

Equation (\ref{sholem}) implies that $x_2$ divides
$\bar R_n^{(1)}(x_2,\ldots,x_n)$; if we set  
$$
\bar R_n^{(1)}(x_1,\ldots,x_{n-1}) = x_1 
R^{(1)\prime}_n(x_1,\ldots,x_{n-1}),
$$ 
where $R^{(1)\prime}_n(x_1,\ldots,x_{n-1})$ belongs to 
$S(\cF\cL_n)$, then 
$R^{(1)\prime}_n(x_1,\ldots,x_{n-1})$ 
is homogeneous of multidegree 
$\delta_2 + \cdots + \delta_{n-1}$.  So 
$R^{(1)\prime}_n(x_1,\ldots,x_{n-1})$ actually belongs to 
$S(\cF\cL_{n-1})$, and may be written  
$S^{(1)}_n(x_2,\ldots,x_{n-1})$.  
We have then 
$$ 
S_n^{(1)}(x_3,\ldots,x_n) x_{n+1} 
= 2 x_3 S_n^{(1)}(x_4,\ldots,x_{n+1}) ; 
$$
this equation is the same as (\ref{sholem}), where
the number of variables is decreased by $1$. Repeating the
reasoning above, we find $\bar R_n^{(1)}(x_1,\ldots,x_{n-1} ) 
= \la x_1\cdots x_{n-1}$, where $\la$ is scalar. Then
equation (\ref{sholem}) implies that $\la = 0$. 
This is a contradiction with $\al >0$. 

Therefore $\al = 0$, which means that $R_n^{(1)}$ is scalar. 
The only cases when a scalar belongs to $FA_{n-1}[[\hbar]]$
is $n = 1$, or this scalar is zero.  
We have therefore shown that if $n>1$, then $R_n^{(1)}$ is zero. 
Equation (\ref{luria}) then implies that $R_n(x_1,\ldots,x_n)$
also vanishes. On the other hand, when $n = 1$, all the solutions of 
(\ref{cordovero}) are $R_n(x) = \la x$, where $\la\in\KK[[\hbar]]$. 
 
Therefore the only solutions to equation (\ref{nasi}) 
are such that if $n>1$, 
$$
\sum_{\al} P^{(n)}_\al(x_1,\ldots,x_n)
\otimes Q^{(n)}_\al(y_1,\ldots,y_n)
$$ is zero, and 
$\sum_\al P^{(1)}_\al(x_1)\otimes Q^{(1)}_\al(y_1)$
is of the form $\la x_1\otimes y_1$, with $\la\in\KK[[\hbar]]$.
This solution corresponds to the assignment
\begin{equation} \label{pobeda}
\A\mapsto \rho_\A = \la\id_\A.
\end{equation} 
So all assignments of 
$(\A\mapsto \rho_\A)$ of $\cG_0(\KK)$ necessarily have the 
form (\ref{pobeda}).  The necessary and
sufficient condition for an assignment of the 
form (\ref{pobeda}) to actually belong to $\cG_0(\KK)$
is that $\la\in(1+\hbar\KK[[\hbar]])$. This ends 
the proof of Theorem \ref{thm:G0}. \hfill \qed\medskip

\begin{remark}
It is much simpler to prove Theorem \ref{thm:G0} by noting that each 
homogeneous component of the right hand side of equation 
(\ref{univ:bar:lev}) is antisymmetric in its two last tensor factors. 
However, the techniques of the above proof will again be used in next 
proofs.  
\end{remark}

\section{Proofs of Theorem \ref{thm:ders} and Proposition \ref{prop:ders:free}}

\subsection{Proof of Theorem \ref{thm:ders}}

Let $(\A\mapsto\la_\A)$ be an element of $\cD$. According to the
`non-$\hbar$-adically completed' version of Proposition
\ref{prop:desc:E}, $(\A\mapsto\la_\A)$ is the image by $i$ of an
element $q$ of $\oplus_{n|n\geq 1} (FL_n\otimes FL_n)_{\SG_n}$.
Let us write $q$ as follows 
$$
q = \sum_{n|n\geq 1} \sum_{\al} P_n^{\prime (\al)}(x_1,\ldots,x_n)
\otimes Q_n^{\prime (\al)}(y_1,\ldots,y_n) ;  
$$
then in the same way as identity (\ref{univ:bar:lev}), one shows
that 
\begin{align} \label{univ:bar:lev'}
& \sum_{n|n\geq 1} \sum_\al  
[P^{\prime (n)}_\al(x_1,\ldots,x_n),x]\otimes 
Q^{\prime (n)}_\al(y_1,\ldots,y_n)  \otimes y
\\ & \nonumber 
- [P^{\prime (n)}_\al(x_1,\ldots,x_n),x]\otimes 
y \otimes Q^{\prime (n)}_\al(y_1,\ldots,y_n) 
\\ & \nonumber 
= \sum_{n|n\geq 1} \sum_\al  
P^{\prime (n)}_\al(x_1,\ldots,x_n) \otimes x \otimes 
[y,Q^{\prime (n)}_\al (y_1,\ldots,y_n)]
\\ & \nonumber
+ P^{\prime (n)}_\al(x_1,\ldots,x_n) \otimes 
[x,Q^{\prime (n)}_\al(y_1,\ldots,y_n)] \otimes y  
\end{align}
holds in $F^{(abb)}$. Applying the Lie bracket to the two last 
tensor factors of this identity, we get 
\begin{align} \label{block} 
& 2 \sum_{n|n\geq 1} \sum_\al
[P^{\prime(n)}_\al(x_1,\ldots,x_n),x]\otimes 
[Q^{\prime(n)}_\al(y_1,\ldots,y_n),y]
\\ & \nonumber 
= \sum_{n|n\geq 1} \sum_{k = 1}^n \sum_\al 
P^{\prime(n)}_\al(x_1,\ldots,[x,x_k],\ldots,x_n) \otimes 
Q^{\prime(n)}_\al (y_1,\ldots,[y,y_k]\ldots,y_n) . 
\end{align}
Let us separate the homogeneous components of this equation, and 
let us denote by $R'_n$ the element of $FL_n$ such that the identity 
$R'_n(x_1,\ldots,x_n)\otimes y_1\cdots y_n =\sum_\al
P'_\al(x_1,\ldots,x_n)\otimes Q'_\al(y_1,\ldots,y_n)$
holds in $(FL_n\otimes FA_n)_{\SG_n}$. Then $R'_n$ satisfies the identity
\begin{equation} \label{new:cordovero}
[R'_n(x_1,\ldots,x_n),x_{n+1}] - [R'_n(x_2,\ldots,x_{n+1}),x_1]
= \sum_{k = 1}^n R'_n(x_1,\ldots,[x_k,x_{k+1}],\ldots,x_{n+1})
\end{equation}
in $FL_{n+1}$.

If $n = 1$, then (\ref{new:cordovero}) implies that $R'_n = 0$. 
If $n = 2$, then the solutions of
(\ref{new:cordovero}) are of the form $R'_n(x,y) =\la[x,y]$, where 
$\la$ is any scalar. 
  
Let us assume that $n>2$. We will show that the only solution to 
(\ref{new:cordovero}) is $R'_n = 0$. 

Let us proceed as above and introduce the elements $R^{\prime (n)}_i$
of $FA_{n-1}$ ($i = 1,\ldots,n$), such that 
$$
R'_n(x_1,\ldots,x_n) = \sum_{i = 1}^n x_i R^{\prime (n)}_i(x_1,\ldots,x_{i-1},
x_{i+1},\ldots,x_n). 
$$
Let us select in (\ref{new:cordovero}) the terms `starting with $x_1$'. 
We obtain 
\begin{align} \label{new:luria}
& R_n^{\prime (1)}(x_2,\ldots,x_n) x_{n+1} + R'_n(x_2,\ldots,x_{n+1})
\\ & \nonumber 
= x_2 R_n^{\prime (1)}(x_3,\ldots,x_{n+1}) 
+ \sum_{i = 2}^n R_n^{\prime (1)}(x_2,\ldots,[x_i,x_{i+1}],\ldots,x_{n+1}). 
\end{align}
There exists unique scalars $\la$ and $(\la_{ij})_{2\leq i<j\leq n}$, 
such that 
\begin{align*}
& R'_n(x_2,\ldots,x_n) = \la x_2\cdots x_n 
+ \sum_{2\leq i<j\leq n} x_2 \cdots x_{i-1}
[x_i,x_j] x_{i+1} \cdots x_{j-1} x_{j+1} \cdots x_{n} 
\\ & +  \on{terms\ of\ degree \ }<n-2 
\end{align*}
(the degree is with respect to the enveloping algebra 
filtration of $\cF\cA_{n-1} = U(\cF\cL_{n-1})$).  We therefore 
obtain the equality 
\begin{align} \label{penn}
& \sum_{i,j|2\leq i<j\leq n} \la_{ij} 
x_2\cdots x_{i-1}[x_i,x_j]x_{i+1} \cdots x_{j-1} x_{j+1}\cdots x_{n+1}
\\ & \nonumber 
= \sum_{i,j|2\leq i<j\leq n} \la_{ij}
x_2\cdots x_{i}[x_{i+1},x_{j+1}]x_{i+2} \cdots x_{j} x_{j+2}\cdots x_{n+1}  
\\ & \nonumber 
+ \la \sum_{i = 2}^n x_2\cdots x_{i-1}[x_i,x_{i+1}]x_{i+2} \cdots x_{n+1} 
\end{align}
modulo terms of degree $\leq \on{max}(1,n-2)$, which is $\leq n-2$ by 
assumption on $n$. 

\begin{lemma}
If $R'_n\neq 0$, then $\la\neq 0$. 
\end{lemma}

{\em Proof.}
Let $d$ be the degree of $R^{\prime(1)}_n$ for the enveloping algebra 
filtration of $\cF\cA_{n-1} = U(\cF\cL_{n-1})$. If $d = 0$, then 
$R^{\prime (1)}_n$ is a scalar, so $n = 1$, which we ruled out. 
So $d>0$.  
Let us denote by $\overline{R^{\prime(n)}_1}$ the image of 
$R^{\prime (n)}_1$ in the degree $d$ part of the associated graded of 
$\cF\cA_{n-1}$. Then since $R'_n(x_2,\ldots,x_{n+1})$ has degree $\leq 1
\leq d$ , and $R^{\prime (1)}_n(x_2,\ldots,[x_i,x_{i+1}],\ldots,x_{n+1})$
has degree $d$, the image of (\ref{new:luria}) in the degree 
$d+1$ part of the associated graded of $\cF\cA_{n-1}$ yields
$$
\overline{R^{\prime(n)}_1}(x_2,\ldots,x_n) x_{n+1}
= x_{2} \overline{R^{\prime(n)}_1}(x_3,\ldots,x_{n+1}).  
$$
Therefore $x_2$ divides 
$\overline{R^{\prime(n)}_1}(x_2,\ldots,x_n)$, and an induction 
as above shows that there exists a scalar $\al$ such that 
$$
\overline{R^{\prime(n)}_1}(x_2,\ldots,x_n) = \al \times\on{\ class\ of\ } 
x_2 \cdots x_n.
$$ 
Since $\overline{R^{\prime(n)}_1}$ cannot be zero, 
$\al$ is not equal to zero. 
On the other hand, we have necessarily $\al = \la$, so $\la\neq 0$. 
\hfill \qed\medskip

Let us assume that $R'_n\neq 0$. We have seen that then $\la\neq 0$.
On the other hand, the image of  (\ref{penn}) in the associated graded of 
$\cF\cA_{n-1}$ implies the equalities 
$$ 
\la_{23} = \la, \quad \la_{34} - \la_{23} = \la, \quad\ldots, \quad
\la_{n-1,n} -
\la_{n-2,n-1} = \la, \quad -\la_{n-1,n} = \la.  
$$ 
Summing up these equalities, we get $(n-1)\la = 0$, so $\la = 0$, a
contradiction. Therefore $R'_n = 0$. It follows that $q$ is
homogeneous of degree $2$, and is therefore proportional to
$[x_1,x_2]\otimes[y_1,y_2]$. This ends the proof of the first 
part of Theorem \ref{thm:ders} on universal derivations. 

Let us prove the statement on universal coderivations.
Assume that $\A\mapsto \la_\A$ is a universal coderivation. 
Then the assignment $\A\mapsto (\la_{\A^*})^t$ is a universal 
derivation. We have shown that a universal derivation is necessarily 
proportional to $u_\A = [,]_\A\circ\delta_\A$. Since
$(u_{\A^*})^t = u_\A$, any universal coderivation is also proportional 
to $\A\mapsto u_\A$. This implies the second part of 
Theorem \ref{thm:ders}.

\subsection{Proof of Proposition \ref{prop:ders:free}}

We have shown that any functorial assignement $(\C\mapsto\la_\C)$, 
where for each object $\C$ of $\LCA$, $\la_\C$ is a derivation of 
$F(\C)$, is provided by an element $q = (q_n)_{n\geq 1}$ of 
$\oplus_{n|n\geq 1}
(FL_n\otimes FL_n)_{\SG_n}$. More precisely, $\la_\C$
is uniquely determined by its restriction to $\C$, which 
has the form (if $\C$ is finite-dimensional) 
$$
\la_\C(x) = \sum_{n\geq 1} \sum_\al 
\sum_{i_1,\ldots,i_n\in I} \langle x, S^{(n)}_\al(b_{i_1},\ldots,b_{i_n})
\rangle R^{(n)}_\al(a_{i_1},\ldots,a_{i_n}), 
$$
where $q_n$ is the class of $\sum_\al R^{(n)}_\al\otimes S^{(n)}_\al$ and 
we write the canonical element of $\C\otimes\C^*$ 
as  $\sum_{i\in I} a_i\otimes b_i$. 

Since $\la_\C$ is a coderivation, it satisfies the 
identity 
\begin{equation} \label{tsedaka}
\delta_{F(\C)} \circ \la_\C(x)
= (\la_\C \otimes \id + \id\otimes \la_\C)\circ\delta_\C(x)  
\end{equation}
for any $x\in\C$. 

Let us denote by $\ad^*$ the coadjoint action of $\C^*$ on $\C$. 
For any $x\in\C$, we have $\delta_\C(x) = \sum_{i\in I}
a_i\otimes \ad^*(b_i)(x)$. 
Then 
\begin{align*} 
& (\la_\C\otimes\id)\circ \delta_\C(x) 
\\ & 
= \sum_{n\geq 1} \sum_\al \sum_{i',i_1,\ldots,i_n\in I} 
\langle a_{i'} , S^{(n)}_\al(b_{i_1},\ldots,b_{i_n})\rangle
R^{(n)}_\al(a_{i_1},\ldots,a_{i_n}) \otimes 
\ad^*(b_{i'})(x)
\\ & 
= \sum_{n\geq 1} \sum_\al \sum_{i',i_1,\ldots,i_n\in I} 
R^{(n)}_\al(a_{i_1},\ldots,a_{i_n})  
\otimes \ad^*\big(S^{(n)}_\al(b_{i_1},\ldots,b_{i_n})\big)(a_{i'}) 
\langle x, b_{i'} \rangle . 
\end{align*} 
On the other hand, if we again denote by $\ad^*$ the action of 
$\C^*$ on $F(\C)$ induced by the coadjoint action of $\C^*$ on $\C$, 
we have 
$$
\delta_{F(\C)} \circ \la_\C(x)  
= \sum_{n\geq 1} \sum_\al \sum_{i',i_1,\ldots,i_n\in I}
a_{i'} \otimes \ad^*(b_{i'})
\big( R^{(n)}_\al(a_{i_1},\ldots,a_{i_n}) \big) 
\langle x, S^{(n)}_\al(b_{i_1},\ldots,b_{i_n})\rangle. 
$$ 
Identity (\ref{tsedaka}) therefore implies 
\begin{align*} 
& \sum_{n\geq 1} \sum_\al \sum_{i',i_1,\ldots,i_n\in I}
a_{i'} \otimes \ad^*(b_{i'})
\big( R^{(n)}_\al(a_{i_1},\ldots,a_{i_n}) \big) 
\otimes S^{(n)}_\al(b_{i_1},\ldots,b_{i_n}) 
\\ & 
= \sum_{n\geq 1} \sum_\al \sum_{i',i_1,\ldots,i_n\in I} 
R^{(n)}_\al(a_{i_1},\ldots,a_{i_n})  
\otimes \ad^*\big(S^{(n)}_\al(b_{i_1},\ldots,b_{i_n})\big)(a_{i'}) 
\otimes b_{i'} 
\\ &   
- \ad^*\big(S^{(n)}_\al(b_{i_1},\ldots,b_{i_n})\big)(a_{i'}) 
\otimes 
R^{(n)}_\al(a_{i_1},\ldots,a_{i_n})  
\otimes b_{i'}  . 
\end{align*}
Let us assume that $\C$ is a Lie bialgebra. Then there
is a unique Lie algebra morphism $\al_\C$ from $F(\C)$ to $\C$, extending the
identity on $\C$. The image of this identity by $\al_\C\otimes\al_\C\otimes
\id$ is the identity in $D(\C)\otimes D(\C)\otimes\C^*$ 
(where both sides belongs to $\C\otimes\C\otimes\C^*$)
\begin{align*} 
& \sum_{n\geq 1} \sum_\al \sum_{i',i_1,\ldots,i_n\in I}
a_{i'} \otimes [b_{i'} , 
R^{(n)}_\al(a_{i_1},\ldots,a_{i_n}) ] 
\otimes S^{(n)}_\al(b_{i_1},\ldots,b_{i_n}) 
\\ & 
+ [a_{i'},R^{(n)}_\al(a_{i_1},\ldots,a_{i_n}) ] 
\otimes b_{i'}  
\otimes S^{(n)}_\al(b_{i_1},\ldots,b_{i_n}) 
\\ & 
= \sum_{n\geq 1} \sum_\al \sum_{i',i_1,\ldots,i_n\in I} 
- R^{(n)}_\al(a_{i_1},\ldots,a_{i_n})  
\otimes a_{i'} \otimes 
[S^{(n)}_\al(b_{i_1},\ldots,b_{i_n}),b_{i'}] 
\\ & 
+ a_{i'} \otimes R^{(n)}_\al(a_{i_1},\ldots,a_{i_n})   
\otimes [S^{(n)}_\al(b_{i_1},\ldots,b_{i_n}), b_{i'}]  ;  
\end{align*}
this identity holds only due to the fact that $\sum_{i\in I} 
a_i\otimes b_i$ is a solution of CYBE, so it holds at the universal 
level. It means that
$q^{(21)}$ satisfies (\ref{univ:bar:lev'}). 
The proof of Theorem \ref{thm:ders} then implies that $q$
is homogeneous of degree $2$, which is the conclusion of
Proposition \ref{prop:ders:free}.

\subsection{Proof of Proposition \ref{isotropy}}

Let $P= (P_n)_{n\geq 1}$ be an element of $\cG(\KK)$, 
such that $P * B = B$. Multiplying $P$ by the suitable power of 
$S_B^2$, we may assume that $P_2 = 0$. The neutral element of 
$\cG(\KK)$ is the sequence $e = (e_n)_{n\geq 1}$, where $e_i = 0$
for $i\geq 2$. Assume that $P$ is not equal to $e$ and let 
$k$ be the smallest index such that $P_k\neq 0$. Then $k\geq 3$. 

Let us denote by $\Delta_0$ the usual (undeformed) coproduct
of $T(\C)[[\hbar]]$, and by $\Delta_1$ the first jet of its deformation: 
so $\Delta_1$ is the unique map from $T(\C)[[\hbar]]$
to $T(\C)^{\otimes 2}[[\hbar]]$, such that $\Delta_{1|\C} = \delta_\C$
and $\Delta_1(xy) = \Delta_0(x)\Delta_1(y) + \Delta_1(x)\Delta_0(y)$
for any pair $x,y$ of elements of $T(\C)[[\hbar]]$. 
 
Then we have the identities 
$$
\Delta_0(\delta_\C^{(P_k)}(x)) = \delta_\C^{(P_k)}(x)\otimes 1 
+ 1\otimes \delta_\C^{(P_k)}(x) , 
$$
and 
$$
\Delta_0(\delta_\C^{(P_{k+1})}(x)) + \Delta_1(\delta_\C^{(P_k)}(x))
= \delta_\C^{(P_{k+1})}(x)\otimes 1 
+ 1\otimes \delta_\C^{(P_{k+1})}(x) + (\delta_\C^{(P_k)} \otimes
\id + \id\otimes \delta_\C^{(P_k)})\circ\delta_\C(x).  
$$
The first identity means that $\delta_\C^{(P_k)}(x)$ is actually
contained in $F(\C)$. Let us expand $P_k$ in the 
form $P_k = \sum_{\sigma\in\SG_k} P_{k,\sigma} x_{\sigma(1)}
\cdots x_{\sigma(k)}$, then this means that   
$P_k = \sum_{\sigma\in\SG_k} P_{k,\sigma^{-1}} x_{\sigma(1)}
\cdots x_{\sigma(k)}$ also belongs to $FL_k[[\hbar]]$. Therefore
the class of 
$\sum_{\sigma\in\SG_k} P_{k,\sigma} x_{\sigma(1)}\cdots x_{\sigma(k)}
\otimes y_1\cdots y_k$ in $(FA_n\otimes FA_n)_{\SG_n}$ actually
belongs to $(FL_n\otimes FL_n)_{\SG_n}$. 

In the second identity, the first and last terms are 
antisymmetric, while the others are symmetric. It follows that 
$$
\Delta_1(\delta_\C^{(P_k)}(x))
= (\delta_\C^{(P_k)} \otimes \id 
+ \id\otimes \delta_\C^{(P_k)})\circ\delta_\C(x).  
$$
Let us denote by $\delta_{F(\C)}$ the extension of $\delta_\C$
to a cocycle map from $F(\C)$ to $\wedge^2 F(\C)$. 
Then the restriction of $\Delta_1$ to $F(\C)$ coincides with 
$\delta_{F(\C)}$, so 
$$
\delta_{F(\C)} \circ \delta_\C^{(P_k)}(x)
= (\delta_\C^{(P_k)} \otimes \id 
+ \id\otimes \delta_\C^{(P_k)})\circ\delta_\C(x).  
$$
for any $x\in\C$. So $\delta_\C^{(P_k)}$ satisfies identity 
(\ref{tsedaka}). Since $k\geq 3$, the proof of Proposition 
\ref{prop:ders:free} implies that $\delta^{(P_k)}_\C = 0$, so 
$P_k = 0$, a contradiction. This proves Proposition \ref{isotropy}.

\subsection*{Acknowledgements}

I would like to thank P.\ Etingof and  C.\ Reutenauer for 
discussions related to this work.

\end{document}